\documentclass{article}
\usepackage{amsmath,amsfonts,amsthm}
\newtheorem{thm}{Theorem}[section]
\newtheorem{pro}[thm]{Proposition}
\newtheorem{cor}[thm]{Corollary}
\newtheorem{lem}[thm]{Lemma}
\newtheorem{conjecture}[thm]{Conjecture}

\newtheorem{rem}[thm]{Remark}
\newtheorem{defn}[thm]{Definition}

\def\div{\raise 1pt \hbox{\big|}}
\def\ei{e_{\infty}}

\begin{document}

\title{Exponents and the Cohomology of Finite Groups}
\author{Jonathan Pakianathan}
\maketitle
\centerline{Dept. of Mathematics}
\centerline{University of Wisconsin}
\centerline{Madison, WI 53705.}

\begin{abstract}
We will provide an example of a $p$-group G which has elements of order
$p^3$ in some of its integral cohomology groups but which also has the property
that $p^2$ annihilates $\bar{H}^i(G;\mathbb{Z})$ for all sufficiently high
$i$. This provides a counterexample to a conjecture of A. Adem which stated
that if a finite group K has an element of order $p^n$ in one of its integral
cohomology groups then it has such an element in infinitely many of its
cohomology groups.

\end{abstract}

\section{Introduction}

Throughout this paper, we will use the integers $\mathbb{Z}$ as 
coefficients
for cohomology groups unless otherwise specified and will write 
$H^*(\cdot)$ for $H^*(\hspace{0.05in} \cdot \hspace{0.05in} ; \mathbb{Z})$.

It is well known that for $G$ a finite group, the integral cohomology groups
$H^*(G)$ are finitely generated in each dimension and are 
annihilated by $|G|$ in positive dimensions. (Here $|G|$ stands for the order
of $G$.)
Thus if we define $\bar{H}(G) = \oplus_{i=1}^{\infty} H^i(G)$,
we have $|G| \cdot \bar{H}(G) = 0 $.

\begin{defn} Given a group $G$, we define the exponent of $G$ as
$
exp(G) = min \{ n \geq 1 : g^n=1, \forall g \in G \}$. We use the convention
that $exp(G) = \infty$ if the set that we are minimizing over is empty.
\end{defn}

\begin{defn} For a finite group $G$, $e(G)$ is defined to be $exp(\bar{H}(G))$.
It follows easily that $e(G) \div |G|$.
\end{defn}

Let $p$ be a prime and $P$ be a finite $p$-group. It is known that the value
of $e(P)$ (which will be a power of $p$) contains information about the
structure of $P$. In particular there is the following theorem of A. Adem:

\begin{thm}[A. Adem]
If $P$ is a finite $p$-group, then $e(P)=p$ if and only if
 $P$ is an elementary abelian $p$-group.
\end{thm}
(Also note that if $P$ is a $p$-group, then $e(P)=1$ if and only if $P=1$. 
See for example page 149 of \cite{Brown}.)

For a general $p$-group $P$, finding $e(P)$ can be quite difficult. Therefore
we define another related quantity which is sometimes easier to calculate.

\begin{defn} For $P$ a finite $p$-group, $\ei (P) = min \{n \geq 1 :
n\bar{H}(P) \text{ is finite} \} .$
\end{defn} 
Notice that $\ei (P)$ will be a power of $p$ and
$\ei (P) \div e(P) \div |P|$.

\begin{rem} It is easy to see that for $C=\mathbb{Z}/p^n\mathbb{Z}$, the
cyclic group of order $p^n$, one has
$exp(C)=\ei (C)=e(C)=|C|=p^n.$
\end{rem}

A question one is lead to ask is, does $\ei (P) = e(P)$? This is part of a
conjecture of A. Adem stated on page 438 of \cite{Lecture}:

\begin{conjecture}
If $S$ is a finite group, and $H^i(S)$ contains elements of
order $p^n$ for some $i$ then it does so for infinitely many $i$.
In particular for a $p$-group $P$, $\ei (P) = e(P)$.
\end{conjecture}

This conjecture is true in certain cases as is seen in the next proposition:

\begin{pro} Let $P$ be a $p$-group.  \\
If $\ei (P) = 1$ then $e(P)=1$ and $P=1$. \\
If $\ei (P) = p$ then $e(P)=p$ and $P$ is elementary abelian.
\end{pro}
\begin{proof}
The first part follows from standard Nakayama-Rim Theory. (see page 140 of
\cite{Brown}.) The second part follows essentially from the theorem of A. Adem
stated before. (See \cite{Adem} or \cite{Le}.)
\end{proof}

However, it turns out the conjecture is false in general. (At least when
$p$ is odd.) 
The main purpose of this paper is to provide a counterexample to the 
conjecture which will be done in the next section. However before doing 
this, let us prove a few basic properties of $\ei (P)$.

\begin{pro} If $P_1 \leq P_2$ where $P_2$ is a $p$-group then
$\ei (P_1) \div \ei (P_2)$.
\label{pro: divide}
\end{pro}

\begin{proof}
By the Evens-Venkov Theorem, $H^*(P_1)$ is a finitely generated
$H^*(P_2)$-module, say with generators $x_1, \dots, x_n \in
H^*(P_1)$. Let $t=max \{ dim(x_i) : 1 \leq i \leq n \} $ and
$s \in \mathbb{N}$ be such that $\ei (P_2) \cdot H^i(P_2) = 0 $
for $i > s$. Then it is easy to see that for $j > s+t$, we have
$\ei (P_2) \cdot H^j(P_1) = 0$. Thus
$\ei (P_2) \cdot \bar{H}(P_1)$ is finite and $\ei (P_1) \div \ei(P_2)$.
\end{proof}

\begin{cor} If $P$ is a finite $p$-group, then 
$$
exp(P) \div \ei (P) \div e(P) \div |P|.
$$
Furthermore, there are examples of $p$-groups which show that these quantities
are different in general.
\label{cor: fundamental}
\end{cor}
\begin{proof}
For the first part, it only remains to show $exp(P) \div \ei (P)$.
Notice there is a cyclic subgroup $C$ of $P$ of size $exp(P)$.
Thus by proposition~\ref{pro: divide}, one has
$exp(P)=|C|=\ei (C) \div \ei(P)$. Thus we have the first part.
When $P$ is elementary abelian of rank greater than one, then $e(P)=p$
is not equal to $|P|$. When $P$ is extraspecial with $exp(P)=p$, then
by the theorem of A. Adem quoted before, we can see $\ei (P) \neq p$
and hence $\ei (P)$ is not equal to $exp(P)$. Finally, the counterexample
in the next section gives an example of a $p$-group $P$ where $\ei (P)$ is 
not equal to $e(P)$. 
\end{proof}

\begin{defn} Let $S(p^n)$ be the Sylow $p$-group of the symmetric group
on $p^n$ letters.
\end{defn}

\begin{rem} $exp(S(p^n))=\ei (S(p^n)) = e(S(p^n)) = p^n $.
\end{rem}
\begin{proof}
The cyclic group $\mathbb{Z}/p^n\mathbb{Z}$ acts faithfully on itself
by left multiplication and hence embeds in the symmetric group on $p^n$
letters. Thus $S(p^n)$ has a cyclic subgroup of order $p^n$ and hence
$p^n \div exp(S(p^n))$. So in light of corollary~\ref{cor: fundamental},
it remains only to show that $e(S(p^n)) \div p^n$. This follows by an
induction. When $n=1$, $S(p^n)$ is cyclic of order $p$, so this case follows.
The induction proceeds by noting that $S(p^n)$ is isomorphic to the wreath 
product of $S(p^{n-1})$ with $\mathbb{Z}/p\mathbb{Z}$. 
Thus $S(p^n)$ has a subgroup of index $p$
which is isomorphic to the direct product of $p$ copies of $S(p^{n-1})$,
and hence $e(S(p^n)) \div pe(S(P^{n-1})) \div pp^{n-1} = p^n$ by an easy 
transfer argument.
\end{proof}

This allows us to prove the following lemma (given in \cite{Le}, we include
a proof for completeness).

\begin{lem}
Let P be a finite $p$-group, if the intersection of all subgroups of P
of index $p^n$ is trivial then $\ei (P) \div p^n$.
\label{lem: last}
\end{lem}
\begin{proof}
Let $H_1,\dots,H_k$ be the distinct subgroups of index $p^n$ in P. 
Then for each $1 \leq i \leq k$, P acts on the left cosets of $H_i$ and this gives us a homomorphism 
$$ 
\phi_i : P \rightarrow S(p^n),
$$ 
with kernel lying in $H_i$.
Putting these homomorphisms together we get a homomorphism from P into
L which is the direct product of $k$ copies of $S(p^n)$ and this is injective as its
kernel is the intersection of the kernels of the $\phi_i$ maps for all
$1 \leq i \leq k$ which is trivial by assumption. Thus $P$ can be considered
a subgroup of $L$ and from the remark above, $\ei (L) = p^n$ and hence
the lemma follows from proposition~\ref{pro: divide}.
\end{proof}

Neither lemma~\ref{lem: last} or proposition~\ref{pro: divide} are true
if we replace $\ei (\cdot)$ with $e(\cdot)$. For example, the group 
$G(\mathfrak{sl}_2)$ 
provided as a counterexample in the next section, embeds in a direct product
$L$ of a few copies of $S(p^2)$ because the intersection of its index $p^2$
subgroups is trivial. Of course $e(L)=p^2$, however 
$e(G(\mathfrak{sl}_2))=p^3$.

\section{The Counterexample}

Throughout this section, $p$ will be an odd prime.
We will be looking at a certain central extension:

$$
1 \rightarrow W \rightarrow G \rightarrow V \rightarrow 1
$$

where $V,W$ are elementary abelian $p$-groups. It is well-known
that there is a bracket $\langle \cdot,\cdot \rangle : V \rightarrow W$ and a 
$p$-power map $\phi : V \rightarrow W$ which are defined using the commutator and
$p$-power in the group $G$ (See \cite{BrP}). 
These are an alternating bilinear form and
a linear map respectively. We will assume as in \cite{BrP} that
$\phi$ is an isomorphism and hence identify $V$ and $W$.
It was shown there that such groups $G$ are in natural correspondence
with bracket algebras (Lie algebras minus the Jacobi identity)
 over $\mathbb{F}_p$. In particular, to every such bracket algebra, there
exists a unique such group corresponding to it. To get the corresponding
bracket algebra one just forms $[\cdot, \cdot] : V \rightarrow
V$ by composing $\langle \cdot, \cdot \rangle $ with the inverse of $\phi$. Hence either $V$ or $W$ can be considered as
the bracket algebra. 

Now recall we have the Lie algebra $\mathfrak{sl}_2$ which is a 3-dimensional
algebra over $\mathbb{F}_p$. We can choose a basis $\{ h, x_+, x_- \}$ which
then has bracket given by:

\begin{align*}
\begin{split}
[h, x_+] &= 2x_+ \\
[h, x_-] &= -2x_- \\
[x_+,x_-] &= h 
\end{split}
\end{align*}

Let $G=G(\mathfrak{sl}_2)$ be the group associated to this Lie algebra as mentioned above. It follows easily then that $G$ has exponent $p^2$ and order
$p^6$. It was shown in the last section of \cite{BrP} that 
$e(G)=p^3$ and in fact there
are elements of order $p^3$ in $H^4(G)$.

We will show that the intersection of all
subgroups of index $p^2$ in $G$ is trivial and hence 
that $\ei (G) \div p^2$ by lemma~\ref{lem: last}. One can conclude
that $\ei(G)=p^2$, as $G$ is not elementary abelian. Since we
know as previously remarked that $e(G)=p^3$, this means $H^*(G)$ is
 annihilated by $p^2$ in all
sufficiently high dimensions but not in all positive dimensions, and hence
gives us the counterexample we sought! 

So it remains to show that the intersection of all subgroups of
index $p^2$ in $G$ is trivial.
Well we have

$$
1 \rightarrow W \rightarrow G \rightarrow V \rightarrow 1
$$
where $W$ and $V$ can both be identified with $\mathfrak{sl}_2$.
Note the preimage of any 1-dimensional subspace of $V$ is a subgroup
of index $p^2$ in $G$ and the intersection of these is in $W$ so 
the intersection of all index $p^2$ subgroups
lies in $W$. To show the intersection of the index $p^2$ subgroups is in fact 
trivial, we need a lemma:

\begin{lem}
Every 2-dimensional sub(Lie)algebra $S$ of $\mathfrak{sl}_2$ considered
as a subset in $V$, has a subgroup $K$ of index $p^2$ in $G$ lying over it.
Furthermore the intersection of $K$ and $W$ is just $\phi(S)$ where $\phi$
is the p-power map mentioned before.
\end{lem}
\begin{proof}
Let $\hat{x}$ and $\hat{y}$ be a basis for the 2-dimensional subalgebra and 
lift them to $x,y \in G$. Let $K$ be the subgroup they generate then $K$ maps
down to the subalgebra under the projection map to $V$ and hence lies
over it. It remains to show $K$ has index $p^2$. Well any element $t$ in $K$
can be written as a product of $x,y,x^{-1},y^{-1}$ in some combination or
other. As the commutator of $x$ and $y$ is central we can move all the $x$'s
to the leftmost and all the $y$'s to the right of them leaving a bunch
of commutators of $x$'s and $y$'s and their inverses on the right. However
if we let $L$ be the image under the $\phi$ p-power map of the subalgebra $S$
then $L$ is generated by $x^p$ and $y^p$ so is in $K$, is central and contains
these commutators. This means that each commutator can be written as a product
of some power of $x^p$ and some power of $y^p$ and as these are central
we can move them to join the other powers of $x$ and $y$ respectively.
The upshot is that any $t \in K$ can be written $t=x^ly^s$ for some $l,s$ integers.
However $x,y$ have order $p^2$ so we see the order of $K$ is at most $p^4$,
but we see easily that it has order at least $p^4$ so $K$ has order $p^4$
and is of index $p^2$. The intersection of $K$ with $W$ is obviously $L$.
\end{proof}

If we view $W$ as $\mathfrak{sl}_2$ then we have already argued that
the intersection of all subgroups of index $p^2$ in $G$ lies in $W$. To show
this intersection is trivial, in light of the lemma above, it suffices to show
that the intersection of 2-dimensional Lie subalgebras in $\mathfrak{sl}_2$
is trivial. We do this next and then will be done.

Note $\{h, x_+ \}$ generates a 2-dimensional subalgebra of $\mathfrak{sl}_2$ and
so does $\{h, x_- \}$. Their intersection is the 1-dimensional subspace
spanned by $h$. So we only need to show that there is a 2-dimensional
subalgebra which does not contain $h$.
Now let $\alpha = 4^{-1}$ in $\mathbb{F}_p$. 
Consider $S$ the subspace generated by
$\{ h+x_+, -\alpha h+x_-\}$. Then we have:

\begin{align*}
\begin{split}
[h+x_+,-\alpha h+x_-] &= 2\alpha x_+ -2x_- + h \\
&= 2\alpha (x_+ + h) -2(x_- -\alpha h) +(1-4\alpha )h \\
&= 2\alpha (h + x_+) -2(-\alpha h + x_-) 
\end{split}
\end{align*}
as $4\alpha =1$.

So we see $S$ is a 2-dimensional sub(Lie)algebra of $\mathfrak{sl}_2$.
However it is easy to see that $h$ is not in $S$. 
Thus we are done as mentioned before.

\end{document}